\documentclass[11pt,a4paper]{amsart}
\usepackage{amsfonts,amscd,amssymb,amsmath,amsthm,hyperref,mathrsfs,xcolor,lscape,longtable}
\usepackage{microtype}
\newcommand{\bb}[1]{\mathbb{#1}}

\newcommand{\C}{\bb{C}}
\newcommand{\Q}{\bb{Q}}

\newcommand{\Z}{\bb{Z}}

\newtheorem{thm}{Theorem}[section]

\theoremstyle{definition}

\numberwithin{equation}{section}
\newcommand{\Sp}{\mathrm{Sp}}

\newcommand{\GL}{\mathrm{GL}}

\begin{document}
\date{\today}
\subjclass[2000]{??, ??}
\title{Arithmeticity of Some Hypergeometric Groups}
\author{Jitendra Bajpai, Sandip Singh and Shashank Vikram Singh}
\address{Max Planck Institute for Mathematics, Bonn, Germany}
\email{jitendra@math.uni-goettingen.de}
\address{Department of Mathematics, Indian Institute of Technology Bombay, Mumbai, India}
\email{sandip@math.iitb.ac.in}
\address{Department of Mathematics, Indian Institute of Technology Bombay, Mumbai, India}
\email{shashankvikrams6@gmail.com}
\subjclass[2010]{Primary: 22E40;  Secondary: 32S40;  33C80}  
\keywords{Hypergeometric group, Monodromy representation, Orthogonal group, Symplectic group}
\begin{abstract}
We show that the hypergeometric groups associated to the pairs of the parameters $\left(0,0,\frac{1}{3}, \frac{2}{3}\right)$, $\left(\frac{1}{2},\frac{1}{2},\frac{1}{4},\frac{3}{4}\right)$; and  $\left(0,\frac{1}{12}, \frac{5}{12},\frac{7}{12},\frac{11}{12}\right), \left(\frac{1}{2},\frac{1}{3},\frac{1}{3},\frac{2}{3},\frac{2}{3}\right)$ are arithmetic.
\end{abstract}
\maketitle


\section{Introduction}\label{sec:intro}The hypergeometric group $\Gamma(\alpha,\beta)$ associated to the pair of the parameters $\alpha=(\alpha_1,\alpha_2,\ldots,\alpha_n)$, $\beta=(\beta_1,\beta_2,\ldots,\beta_n)\in\C^n$ (with the condition that $\alpha_j-\beta_k\notin\Z, \forall 1\le j,k\le n$) is defined (up to conjugations in $\GL_n(\C)$) as the subgroup of $\GL_n(\C)$ generated by the companion matrices $A$ and $B$ of the polynomials $f(x)=\prod_{j=1}^n(x-e^{2\pi i\alpha_j})$ and $g(x)=\prod_{j=1}^n(x-e^{2\pi i\beta_j})$, respectively (cf. \cite[Theorem 3.5]{BH}). That is, if we write \[f(x)=x^n+a_{n-1}x^{n-1}+\cdots+a_1x+a_0,\qquad g(x)=x^n+b_{n-1}x^{n-1}+\cdots+b_1x+b_0\]
then

 \[A=\begin{pmatrix}
0&0&\cdots&0&-a_0\\1&0&\cdots&0&-a_1\\0&1&\cdots&0&-a_2\\ \vdots&\vdots&\ddots&\vdots&\vdots\\ 0&0&\cdots&1&-a_{n-1}
\end{pmatrix},\qquad B=\begin{pmatrix}
0&0&\cdots&0&-b_0\\1&0&\cdots&0&-b_1\\0&1&\cdots&0&-b_2\\ \vdots&\vdots&\ddots&\vdots&\vdots\\ 0&0&\cdots&1&-b_{n-1}
\end{pmatrix}\]
and \[\Gamma(\alpha,\beta)=<A,B>\subseteq\GL_n(\C).\]

From here onward we  denote the hypergeometric group corresponding to the pair of the parameters $\alpha, \beta\in\C^n$ by the notation $\Gamma(f,g)$ instead of $\Gamma(\alpha,\beta)$ where the polynomials $f,g$, corresponding to the pair of the parameters $\alpha,\beta$, are defined as above.

We consider the cases where $f,g$ are products of the cyclotomic polynomials. In these cases $f(0)=\pm 1$, $g(0)=\pm 1$, and hence the corresponding companion matrices $A,B\in\GL_n(\Z)$ and $\Gamma(f,g)\subseteq\GL_n(\Z)$. If $f,g$ are also coprime (this condition is required as $\alpha_j-\beta_k\notin\Z,\ \forall 1\le j,k\le n$) and form a primitive pair (that is, $\nexists f_1, g_1\in\Z[x]$ so that $f(x)=f_1(x^k),\ g(x)=g_1(x^k)$ for some $k\ge 2$), then it follows from Beukers and Heckman \cite[Theorem 6.5]{BH} that if $f(0)=g(0)=1$ (and this happens only when $n$ is even) then the corresponding hypergeometric group $\Gamma(f,g)$ preserves a non-degenerate symplectic form $\Omega$ on $\Z^n$ and $\Gamma(f,g)$ is a Zariski dense subgroup of the corresponding symplectic group $\Sp_\Omega$. It also follows from Beukers and Heckman \cite[Theorem 6.5]{BH} that if $\frac{f(0)}{g(0)}=-1$ and $\Gamma(f,g)$ is infinite, then $\Gamma(f,g)$ preserves a non-degenerate quadratic form $\mathrm{Q}$ on $\Z^n$ and $\Gamma(f,g)$ is a Zariski dense subgroup of the corresponding orthogonal group $\mathrm{O}_\mathrm{Q}$.

Note that  if we denote the Zariski closure of $\Gamma(f,g)$ by $\mathrm{G}$, then $\Gamma(f,g)\subseteq \mathrm{G}(\Z)$ (for $f,g$ being products of cyclotomic polynomials) in both cases, and so we call $\Gamma(f,g)$ arithmetic if the index of $\Gamma(f,g)$ inside $\mathrm{G}(\Z)$ is finite, and thin otherwise.

Singh and Venkataramana \cite[Tables 1,2]{SV} provide a list  of $111$ pairs of degree four polynomials $f,g$, satisfying the above conditions. It follows from \cite{SS, SS17, SV} that the hypergeometric groups corresponding to 87 pairs (cf. \cite[Tables 1,3]{SS17}, \cite[Table 1]{SV}) of them are arithmetic and the thinness of the hypergeometric groups corresponding to 13 (cf. \cite[Table 2]{SS17}) of the remaining 24 pairs follows from Brav and Thomas \cite{BT}. 

So far, the question to determine the arithmeticity or thinness of the hypergeometric groups corresponding to the remaining 11 pairs (cf. \cite[Table 4]{SS17}) is open and we know that to answer this question it is enough to show the arithmeticity or thinness of just 6 groups as the other 5 are scalar shifts of 5 of these 6 groups and their arithmeticity or thinness follows (cf. \cite[Remarks 1,2]{SS17}). In this article, we show that one of the remaining 6 independent pairs correspond to an arithmetic hypergeometric group (and as said above  it shows the arithmeticity of the hypergeometric groups corresponding to 2 of the remaining 11 pairs). In fact, we obtain the following theorem.

\begin{thm}\label{firsttheorem}
 The hypergeometric group associated to the pair of parameters $\left(0,0,\frac{1}{3}, \frac{2}{3}\right),$ $\left(\frac{1}{2},\frac{1}{2},\frac{1}{4},\frac{3}{4}\right)$ is arithmetic in the corresponding symplectic group. Since the pair of the parameters $\left(0,0,\frac{1}{4}, \frac{3}{4}\right),$ $\left(\frac{1}{2},\frac{1}{2},\frac{1}{6},\frac{5}{6}\right)$ is a scalar shift of  $\left(0,0,\frac{1}{3}, \frac{2}{3}\right),$ $\left(\frac{1}{2},\frac{1}{2},\frac{1}{4},\frac{3}{4}\right)$, the arithmeticity of the corresponding hypergeometric group follows.  These are Examples 1 and 7 of \cite[Table 4]{SS17}.
\end{thm}

Note that Detinko, Flannery and Hulpke \cite{DFH} also show the arithmeticity of these two groups (see Examples 126 and 141 of  \cite[Table 1]{DFH}) by using the coset enumeration algorithm (implemented in GAP) for computing the index of  $\Gamma(f,g)$ inside $\Sp_4(\Z)$, while we show the arithmeticity of these two groups by showing that they intersect the unipotent groups corresponding to the highest and the second highest roots of the symplectic group $\Sp_4$ non-trivially and using the criterion of Venkataramana \cite[Theorem 3.5]{Ve}.

Bajpai and Singh \cite[Tables 1-7] {BS} provide a list of $77$ (up to scalar shifts) possible pairs of degree five polynomials that are products of cyclotomic polynomials and satisfying the above conditions of Beukers and Heckman \cite[Theorem 3.5]{BH} so that the Zariski closures of the corresponding hypergeometric groups are either finite or the orthogonal groups of non-degenerate quadratic forms of signature $(p,q)$ with $p,q\ge 1$. It follows from the criterion of Beukers and Heckman  \cite[Theorem 3.5]{BH} that only 4 (cf. \cite[Table 5]{BS}) of these pairs correspond to finite hypergeometric groups and only 17 pairs (cf. \cite[Tables 1 and 7]{BS})  correspond to the hypergeometric groups having the Zariski closures of real rank one. The remaining 56  pairs (cf. \cite[Tables 2, 3, 4 and  6]{BS}) correspond to the hypergeometric groups having the Zariski closures of real rank two. It follows from Venkataramana \cite{Ve3} (cf. \cite[Table 2]{BS}) that the hypergeometric groups corresponding to 11 of these 56 pairs are arithmetic and the arithmeticity of 2 other cases follows from Singh \cite{SS0}. Then, Bajpai and Singh \cite[Theorem 2; cf. Table 4]{BS} show that the hypergeometric groups corresponding to  23 of the remaining 43 pairs (corresponding to which the Zariski closures of the hypergeometric groups have real rank two) are arithmetic. In this article we also show that one of the remaining 20 pairs (cf. \cite[Table 6]{BS}) correspond to an arithmetic hypergeometric group. In fact, we obtain the following theorem.

\begin{thm}\label{secondtheorem}
 The hypergeometric group associated to the pair of parameters $\left(0,\frac{1}{12}, \frac{5}{12},\frac{7}{12},\frac{11}{12}\right)$, $\left(\frac{1}{2},\frac{1}{3},\frac{1}{3},\frac{2}{3},\frac{2}{3}\right)$ is arithmetic in the corresponding orthogonal group. This is Example 67 of \cite[Table 6]{BS}.
\end{thm}

We prove this theorem by showing that the corresponding hypergeometric group intersects the unipotent groups corresponding to the highest and the second highest roots of the orthogonal group $\mathrm{O}_5$ non-trivially (cf. \cite[Theorem 3.5]{Ve}).

So in this paper we prove the arithmeticity of the hypergeometric groups associated to Examples 1 and 7 of \cite[Table 4]{SS17} and Example 67 of  \cite[Table 6]{BS}.

\section{Proof of Theorem \ref{firsttheorem}}
It is a hypergeometric group in $\mathrm{Sp_4}$ associated to the  parameters
\begin{center}
	$\alpha=(0,0,\frac{1}{3},\frac{2}{3})$, \;  \; $\beta=(\frac{1}{2},\frac{1}{2},\frac{1}{4},\frac{3}{4})$.
\end{center} 

This is Example 1 of \cite[Table 4]{SS17}. In this case,
\begin{center}
	$f(x)= x^4-x^3-x+1$, $\quad$   $g(x)= x^4+2x^3+2x^2+2x+1$
\end{center}
and $f(x)-g(x)=-3x^3-2x^2-3x$. The corresponding hypergeometric group $\Gamma(f,g)$ is generated by $A$ and $B$, that are, respectively, the companion matrices of the polynomials $f(x)$ and $g(x)$. Let $C=A^{-1}B$. So we get
\[
A =
\begin{pmatrix}
0 & 0 & 0 & -1 \\
1 & 0 & 0 & 1\\
0 & 1 & 0 & 0\\
0 & 0 & 1 & 1\\
\end{pmatrix},\qquad
B =
\begin{pmatrix}
0 & 0 & 0 & -1 \\
1 & 0 & 0 & -2\\
0 & 1 & 0 & -2\\
0 & 0 & 1 & -2
\end{pmatrix},\qquad 
C=A^{-1}B =
\begin{pmatrix}
1 & 0 & 0 & -3 \\
0 & 1 & 0 &  -2\\
0 & 0 & 1 &  -3\\
0 & 0 & 0 & 1\\
\end{pmatrix}.
\] 
By a simple computation we find the symplectic form $\Omega$ (determined upto scalar multiples)  preserved by the corresponding hypergeometric group $\Gamma(f,g)$; the matrix form of $\Omega$ with respect to the standard basis of $\Q^4$ over $\Q$ is
\[
\Omega=\begin{pmatrix}
0 & 1 & -{2}/{3} & {1}/{3} \\
-1 & 0 & 1 &  -{2}/{3}\\
{2}/{3} & -1 & 0 &  1\\
-{1}/{3} & {2}/{3} & -1 & 0\\
\end{pmatrix}
\]
and $\Gamma(f,g) \subseteq \Sp_4(\Omega)$ as a Zariski dense subgroup.

By doing further computations we find that if we consider the matrix
\[
P=
\begin{pmatrix}
1 &-3 & 4 & {8}/{3}\\
{-3}/{7} &-2 & 23 & 0\\
0 &-3 & 22 & 1\\
0 &0 & -1 & 0\\
\end{pmatrix}\]
as the change of basis matrix, then the matrix form of $\Omega$, with respect to the new basis of $\Q^4$ over $\Q$, is \[\qquad
\Omega'=P^t\Omega P=
\begin{pmatrix}
0 & 0 & 0 & {1}/{21}\\
0 & 0 & {8}/{3} & 0\\
0 & {-8}/{3} & 0 & 0\\
{-1}/{21} & 0 & 0 & 0\\
\end{pmatrix}
\]
where $P^t$ denotes the transpose of the matrix $P$.

It can be checked easily that for the symplectic form $\Omega'$ the group of diagonal matrices
 
\[T=\left\{\begin{pmatrix}
t_1&0&0&0\\0&t_2&0&0\\
0&0&t_2^{-1}&0\\
0&0&0&t_1^{-1}
\end{pmatrix}: t_1\neq 0\neq t_2\right\}\]form a maximal torus in $\Sp_4(\Omega')$ and if we denote by ${\chi_i}$ the character of $T$ defined by
 
\[\begin{pmatrix}
t_1&0&0&0\\0&t_2&0&0\\
0&0&t_2^{-1}&0\\
0&0&0&t_1^{-1}
\end{pmatrix}\mapsto t_i, \quad \mbox{for }i=1,2,\]then the set $\Phi=\{\chi_1^2, \chi_1\chi_2, \chi_1\chi_2^{-1}, \chi_2^{2},\chi_1^{-2}, \chi_1^{-1}\chi_2^{-1}, \chi_1^{-1}\chi_2, \chi_2^{-2}\}$ forms the root system of $\Sp_4(\Omega')$ and if we fix the set $\Delta=\{\chi_1\chi_2^{-1},\chi_2^2\}$ of simple roots, then $\chi_1^2$ and $\chi_1\chi_2$ are, respectively, the highest and second highest roots of $\Sp_4(\Omega')$ and the corresponding unipotent groups are

\[U_{\chi_1^2}=\left\{\begin{pmatrix}
1&0&0&c\\0&1&0&0\\
0&0&1&0\\
0&0&0&1
\end{pmatrix}: c\in\Q\right\},\qquad U_{\chi_1\chi_2}=\left\{\begin{pmatrix}
1&0&c&0\\0&1&0&\lambda c\\
0&0&1&0\\
0&0&0&1
\end{pmatrix}: c\in\Q\right\}\]for some non-zero scalar $\lambda$.

Now we show that the Zariski dense subgroup $P^{-1}\Gamma(f,g)P$ of $\Sp_4(\Omega')$ contains some non-trivial elements of the unipotent groups $U_{\chi_1^2}$ and $U_{\chi_1\chi_2}$ which shows the arithmeticity of $P^{-1}\Gamma(f,g)P$ inside $\Sp_4(\Omega')$ using the criterion \cite[Theorem 3.5]{Ve} of Venkataramana; and the arithmeticity of $\Gamma(f,g)$ inside $\Sp_4(\Omega)$ follows.

Let $G$ denote the element $B^3A^2B^2AB^2$ of the hypergeometric group $\Gamma(f,g)$. Then the elements $E_1=P^{-1}CP, E_2=P^{-1}GCG^{-1}P$ and $E_3=P^{-1}A^{-1}CAP$ belong to the group $P^{-1}\Gamma(f,g)P$, for which, the corresponding matrices are

\[
E_1=
\begin{pmatrix}
1 & 0 & 0 & 0 \\
0 & 1 & -1 &  0\\
0 & 0 & 1 &  0\\
0 & 0 & 0 & 1\\
\end{pmatrix},\;
E_2=
\begin{pmatrix}
1 & 0 & 0 & 0 \\
0 & 1 & 0 &  0\\
0 & 1 & 1 &  0\\
0 & 0 & 0 & 1\\
\end{pmatrix},\;
E_3=
\begin{pmatrix}
1 & 168 & -1176 & -56 \\
0 & 64 & -441 &  -21\\
0 & 9 & -62 &  -3\\
0 & 0 & 0 & 1\\
\end{pmatrix}.
\]
A computation shows that if $E_4=E_1^7E_3E_1^{-7}$, $E_5=E_2^{-9}E_4$ and  $E_6=E_1E_5E1^{-1}$, then 

\[
E_7=E_5E_6E_5^{-1}E_6^{-1}=\begin{pmatrix}
1 & 0 & 0 & 1008 \\
0 & 1 & 0 &  0\\
0 & 0 & 1 &  0\\
0 & 0 & 0 & 0\\
\end{pmatrix}\in P^{-1}\Gamma(f,g)P
\]
is a non-trivial element of the unipotent group $U_{\chi_1^2}$, corresponding to the highest root $\chi_1^2$ of $\Sp_4(\Omega')$. Also, if $E_8=E_5^{18}E_7$, then

\[
E_9=E_7^{161}E_8E_6^{-18}=\begin{pmatrix}
1 & 0 & -3024 & 0 \\
0 & 1 & 0 &  -54\\
0 & 0 & 1 &  0\\
0 & 0 & 0 & 0\\
\end{pmatrix}\in  P^{-1}\Gamma(f,g)P
\]
is a non-trivial element of the unipotent group $U_{\chi_1\chi_2}$, corresponding to the second highest root $\chi_1\chi_2$ of $\Sp_4(\Omega')$.

The existence of the non-trivial unipotent elements $E_8$ and $E_9$ in $P^{-1}\Gamma(f,g)P$ corresponding to the highest and the second highest roots of $\Sp_4(\Omega')$ completes the proof of Theorem \ref{firsttheorem}. \qed

\section{Proof of Theorem \ref{secondtheorem}}
Observe that the pair of the parameters of Theorem \ref{secondtheorem} is a scalar shift (by $\frac{1}{2}$) of the pair of the parameters $\alpha=(\frac{1}{2},\frac{1}{12},\frac{5}{12},\frac{7}{12},\frac{11}{12})$,
	$\beta=(0,\frac{1}{6},\frac{1}{6},\frac{5}{6},\frac{5}{6})$. Therefore to prove Theorem \ref{secondtheorem} it is enough to show the arithmeticity of the hypergeometric group corresponding to the pair of parameters $\alpha=(\frac{1}{2},\frac{1}{12},\frac{5}{12},\frac{7}{12},\frac{11}{12}),
	\beta=(0,\frac{1}{6},\frac{1}{6},\frac{5}{6},\frac{5}{6})$. This is what we do here.

In this case
\[f(x)= x^5+x^4-x^3-x^2+x+1, \quad g(x)= x^5-3x^4+5x^3-5x^2+3x-1.\]

The corresponding hypergeometric group $\Gamma(f,g)$ is generated by the matrices $A$ and $B$, that are, respectively, the companion matrices of the polynomials $f(x)$ and $g(x)$. Let $C=A^{-1}B$. Then
\[
A =
\begin{pmatrix}
0 & 0 & 0& 0 & -1 \\
1 & 0 & 0& 0 & -1\\
0 & 1 & 0& 0 & 1\\
0 & 0 & 1& 0 & 1\\
0 & 0 & 0& 1 & -1\\
\end{pmatrix},\qquad
B=
\begin{pmatrix}
0 & 0 & 0& 0 & 1 \\
1 & 0 & 0& 0 & -3\\
0 & 1 & 0& 0 & 5\\
0 & 0 & 1& 0 & -5\\
0 & 0 & 0& 1 & 3\\
\end{pmatrix} ,\qquad 
C=A^{-1}B =
\begin{pmatrix}
1 & 0 & 0 &0 & -4
\\
0 & 1 & 0 &0 &  6\\
0 & 0 & 1 &0 &  -4\\
0 & 0 & 0 &1 & 2 \\
0 & 0 & 0 &0 & -1
\end{pmatrix}.
\] 

Let $v=(C-I)(e_5)$, where $e_5=(0,0,0,0,1)$ and $I$ is the $5\times 5$ identity matrix. Then $v=(-4,6,-4,2,-2)$ and the vectors of the set $\{v,Bv,B^2v,B^3v,B^4v\}$ are linearly independent by \cite[Lemma 2.1]{SS0}. Now, by using the method of \cite{SS0} we find that the matrix form of the quadratic form $Q$ (preserved by the hypergeometric group $\Gamma(f,g)$), with respect to the basis $\{v,Bv,B^2v,B^3v,B^4v\}$ of $\Q^5$ over $\Q$ is

\[q=\begin{pmatrix}
-1 & -2 & -3& -1 & 3 \\
-2 & -1 & -2& -3 & -1\\
-3 & -2 & -1& -2 & -3\\
-1 & -3 & -2& -1 & -2\\
3& -1 & -3& -2 & -1
\end{pmatrix}\]
and a simple computation shows that the matrix form of the quadratic form $Q$, with respect to the standard basis $\{e_1,e_2,e_3,e_4,e_5\}$
of $\Q^5$ over $\Q$, is

\[Q=\begin{pmatrix}
{-19}/{9} & {-17}/{9} & {-10}/{9}& {1}/{9} & {8}/{9} \\
{-17}/{9} & {-19}/{9} & {-17}/{9}& {-10}/{9} & {1}/{9}\\
{-10}/{9} & {-17}/{9} & {-19}/{9}& {-17}/{9} & {-10}/{9}\\
{1}/{9} & {-10}/{9} & {-17}/{9}& {-19}/{9} & {-17}/{9}\\
{8}/{9}& {1}/{9} & {-10}/{9}& {-17}/{9} & {-19}/{9}
\end{pmatrix}\]
It follows that $A^tQA=B^tQB=Q$ where $A^t$ and $B^t$ denote the transposes of the corresponding matrices. Hence the hypergeometric group $\Gamma(f,g)$ preserves the quadratic form $Q$ (determined upto scalar multiples) and $\Gamma(f,g) \subseteq \mathrm{O}_5(Q)$ as a Zariski dense subgroup by \cite{BH}. 

By doing further computations we find that if we consider the matrix
\[
P=
\begin{pmatrix}
1 & 0 & 0 & 0 & {-1}/{2}\\
-4 &-1 & 0 & -3 & -1\\
8 &3 & -1 & 4&2\\
-5 &-2 & 1 & 2 & {1}/{2}\\
0 & 0 & -1 & -5 & -2
\end{pmatrix}\]
as the change of basis matrix, then the matrix form of $Q$, with respect to the new basis of $\Q^5$ over $\Q$, is \[\qquad
Q'=P^tQ P=\begin{pmatrix}
0 & 0 & 0 & 0 & 1\\
0 & 0 & 0 & -1 & 0\\
0 & 0 & -1 & 0 & 0\\
0 & -1 & 0 & 0 & 0\\
1 & 0 & 0 & 0 & 0\\
\end{pmatrix}
\]
where $P^t$ denotes the transpose of the matrix $P$.

It can be checked easily that for the quadratic form $Q'$ the group of diagonal matrices
 
\[T=\left\{\begin{pmatrix}
t_1&0&0&0&0\\0&t_2&0&0&0\\
0&0&1&0&0\\
0&0&0&t_2^{-1}&0\\
0&0&0&0&t_1^{-1}
\end{pmatrix}: t_1\neq 0\neq t_2\right\}\]form a maximal torus in $\mathrm{O}_5(Q')$ and if we denote by ${\chi_i}$ the character of $T$ defined by
 
\[\begin{pmatrix}
t_1&0&0&0&0\\0&t_2&0&0&0\\
0&0&1&0&0\\
0&0&0&t_2^{-1}&0\\
0&0&0&0&t_1^{-1}
\end{pmatrix}\mapsto t_i, \quad \mbox{for }i=1,2,\]then the set $\Phi=\{\chi_1,\chi_2, \chi_1\chi_2, \chi_1\chi_2^{-1}, \chi_1^{-1},\chi_2^{-1}, \chi_1^{-1}\chi_2^{-1}, \chi_2\chi_1^{-1}\}$ forms the root system of $\mathrm{O}_5(Q')$ and if we fix the set $\Delta=\{\chi_2,\chi_1\chi_2^{-1}\}$ of simple roots, then $\chi_1\chi_2$ and $\chi_1$ are, respectively, the highest and second highest roots of $\mathrm{O}_5(Q')$ and the corresponding unipotent groups are

\[U_{\chi_1\chi_2}=\left\{\begin{pmatrix}
1&0&0&c&0\\0&1&0&0&c\\
0&0&1&0&0\\
0&0&0&1&0\\
0&0&0&0&1
\end{pmatrix}: c\in\Q\right\},\qquad U_{\chi_1}=\left\{\begin{pmatrix}
1&0&c&0&\frac{1}{2}c^2\\0&1&0&0&0\\
0&0&1&0&c\\
0&0&0&1&0\\
0&0&0&0&1
\end{pmatrix}: c\in\Q\right\}.\]
Now we show that the Zariski dense subgroup $P^{-1}\Gamma(f,g)P$ of $\mathrm{O}_5(Q')$ contains some non-trivial elements of the unipotent groups $U_{\chi_1\chi_2}$ and $U_{\chi_1}$ which shows the arithmeticity of $P^{-1}\Gamma(f,g)P$ inside $\mathrm{O}_5(Q')$ using the criterion \cite[Theorem 3.5]{Ve} of Venkataramana. This shows the arithmeticity of $\Gamma(f,g)$ inside $\mathrm{O}_5(Q)$.

Let $a=P^{-1}AP$ and $b=P^{-1}BP$. Then
\[
a =
\begin{pmatrix}
-2 & -1/2 & -1& -5 & -3/2 \\
2 & 0 & 4& 14 & 4\\
-2 & -1 & 1& -2 & -1\\
3 & 1 & 1& 7 & 5/2\\
-4 & -1 & -4& -20 & -7\\
\end{pmatrix},\qquad
b=
\begin{pmatrix}
-2 & -1/2 & -2& -10 & -7/2 \\
2 & 0 & 4& 14 & 4\\
-2 & -1 & 1& -2 & -1\\
3 & 1 & 1& 7 & 5/2\\
-4 & -1 & -2& -10 & -3\\
\end{pmatrix} .
\]

Let $c=ba^{-1}$, $E_1=b^{-1}aba^{-1}$ and $E_2=E_1c$. Then 
\[
E_3=E_2b^3E_2^{-1}b^{-3}=
\begin{pmatrix}
1 & 0 & 2 &10 & 2\\
0 & 1 & 0 &0 &  10\\
0 & 0 & 1 &0 &  2\\
0 & 0 & 0 &1 & 0 \\
0 & 0 & 0 &0 & 1
\end{pmatrix}
\]
is a non-trivial unipotent element in $P^{-1}\Gamma(f,g)P$.

Now, let $r=ba^{-6}b^{-1}$, $E_4=cE_3c$, $E_5=crc$ and $E_6=E_4E_5E_4^{-1}E_5^{-1}$. Then it follows that
\[
E_7=cE_6c=
\begin{pmatrix}
1 & 72 & -36 &0& 648\\
0 & 1 & 0 &0 &  0\\
0 & 0 & 1 &0 &  -36\\
0 & 0 & 0 &1 & 72\\
0 & 0 & 0 &0 & 1
\end{pmatrix},\qquad
E_8=E_7E_3^{18}=
\begin{pmatrix}
1 & 72 & 0 &180 & 12960\\
0 & 1 & 0 &0 &  180\\
0 & 0 & 1 &0 &  0\\
0 & 0 & 0 &1 & 72\\
0 & 0 & 0 &0 & 1
\end{pmatrix}
\]
are some more unipotent elements inside the group $P^{-1}\Gamma(f,g)P$.

Let $c_1=a^{-1}b$, $d=c_1b^6c_1$, $p=ab^{-6}a^{-1}$, $E_9=bE_2b^{-1}E_2^{-1}c$, $E_{10}=cdcp$, $E_{11}=cE_{10}c$ and $E_{12}=cE_9cE_3^{-2}$. Then
\[
E_{13}= E_{11}E_{12}E_{11}E_{12}=
\begin{pmatrix}
1 & 0 & 0 &96 & 0\\
0 & 1 & 0 &0 &  96\\
0 & 0 & 1 &0 &  0\\
0 & 0 & 0 &1 & 0 \\
0 & 0 & 0 &0 & 1
\end{pmatrix}\in  P^{-1}\Gamma(f,g)P
\]
is a non-trivial element of the unipotent group $U_{\chi_1\chi_2}$, corresponding to the highest root $\chi_1\chi_2$ of $\mathrm{O}_5(Q')$.

Further, let $E_{14}=E_8E_{13}^{-2}E_3$ and $E_{15}=E_{14}^{48}E_{13}$. Then
\[
E_{16}=E_7^{-48}E_{15}=
\begin{pmatrix}
1 & 0 & 1824 &0 & 1663488\\
0 & 1 & 0 &0 &  0\\
0 & 0 & 1 &0 &  1824\\
0 & 0 & 0 &1 & 0 \\
0 & 0 & 0 &0 & 1
\end{pmatrix}\in  P^{-1}\Gamma(f,g)P
\]
is a non-trivial element of the unipotent group $U_{\chi_1}$, corresponding to the second highest root $\chi_1$ of $\mathrm{O}_5(Q')$.

The existence of the non-trivial unipotent elements $E_{13}$ and $E_{16}$ in $P^{-1}\Gamma(f,g)P$ corresponding to the highest and the second highest roots of $\mathrm{O}_5(Q')$ completes the proof of Theorem \ref{secondtheorem}. \qed

\section*{Acknowledgements}
The work of the second author is supported by the DST-INSPIRE Faculty Fellowship no. DST/INSPIRE/04/2015/000794 and the MATRICS grant no. MTR/2021/000368. This paper is part of the Ph.D. thesis of the third author who would like to thank IIT Bombay for providing support.

\nocite{}

\bibliography{Sp6}
\end{document}